\documentclass[journal]{IEEEtran}
\usepackage[latin9]{inputenc}
\usepackage{multirow}
\usepackage{amsmath}
\usepackage{amssymb}
\usepackage{graphicx}

\makeatletter

\providecommand{\tabularnewline}{\\}

\usepackage{amsfonts}
\ifCLASSINFOpdf
\else
\fi
\makeatother

\begin{document}
\title{State Compensation Linearization and Control}
\author{Quan Quan, Jinrui Ren\thanks{Q. Quan is with the School of Automation Science and Electrical Engineering,
Beihang University, Beijing 100191, China. (e-mail: qq\_buaa@buaa.edu.cn.}\thanks{J. Ren is with the School of Automation, Xi'an University of Posts
\& Telecommunications, Shaanxi, Xi'an 710121, China. (e-mail: renjinrui@xupt.edu.cn).}}
\maketitle
\begin{abstract}
The linearization method builds a bridge from mature methods for linear
systems to nonlinear systems and has been widely used in various areas.
There are currently two main linearization methods: Jacobian linearization
and feedback linearization. However, the Jacobian linearization method
has approximate and local properties, and the feedback linearization
method has a singularity problem and loses the physical meaning of
the obtained states. Thus, as a kind of complementation, a new linearization
method named state compensation linearization is proposed in the paper.
Their differences, advantages, and disadvantages are discussed in
detail. Based on the state compensation linearization, a state-compensation-linearization-based
control framework is proposed for a class of nonlinear systems. Under
the new framework, the original problem can be simplified. The framework
also allows different control methods, especially those only applicable
to linear systems, to be incorporated. Three illustrative examples
are also given to show the process and effectiveness of the proposed
linearization method and control framework. 
\end{abstract}


\section{Introduction}

Linearization, aiming to extract a linear relation from a nonlinear
dynamical system, is a common process in control designs. This is
on one hand caused by the fact that most of the dynamical systems
are nonlinear in nature, and is incentivized by the rich results of
linear control theories and methodologies. On the other hand, humans
are more used to dealing with linear systems especially when deploying
and tuning controllers on real systems. The internal model hypothesis
suggests that the central nervous system constructs models of the
body and the physical world, and that humans use approximate feedforward
plant inversion to interact with linear systems \cite{wolpert1998internal},
whereas experiments in \cite{seigler2018on} showed that such a linear
plant inversion in humans' brains may be less precise when directly
dealing with nonlinear systems. Therefore, linearization is important
in both systematic control design \cite{khalil2002nonlinear,slotine1991applied,jakubczyk1980linearization}
and facilitating human-machine interaction. In the control community,
linear systems are regarded as a class of standard and fundamental
systems, where theories have been nearly fully developed. Although
more advanced nonlinear control methods have been developed rapidly
in recent years, compared with linear methods, they are more difficult
to apply in practice. Linearization thus allows mature results of
linear systems, for example, frequency-domain methods \cite{chen1970introduction},
to be applied for nonlinear systems. In engineering practice, linearization
can be considered as an intermediate step, and the subsequent analysis
and design in the linear domain will become much easier to implement
and understand \cite{chi2019adjacent,gan1999linearization}.

There are two main linearization methods: \emph{Jacobian linearization
(i.e., linearization by Taylor's expansion)} and \emph{feedback linearization
}\cite{khalil2002nonlinear},\cite{slotine1991applied},\cite{8865273}.
Jacobian linearization is the most popular one because it is simple,
intuitive, and easy-to-use.\emph{ }However,\emph{ }the Jacobian linearization
method is an approximate linearization method, which neglects higher-order
nonlinearity terms in Taylor's expansion directly. Moreover, Taylor\textquoteright s
expansion holds locally, and thus the Jacobian linearization has a
local property. Different from the Jacobian linearization, feedback
linearization is an accurate linearization method, which transforms
a nonlinear system into an equivalent linear system of the Brunowski
form, namely multiple chains of integrators. Nevertheless, feedback
linearization is somewhat difficult to apply directly in some situations
\cite{slotine1991applied}. Full state measurement is also required
for implementing feedback linearization, and robustness may not be
guaranteed when there exist uncertain parameters or unmodelled dynamics
\cite{khalil2002nonlinear}. Moreover, feedback linearization may
lead to a singularity problem, i.e., a vanishing denominator, which
makes the system uncontrollable \cite[pp.213-216]{slotine1991applied}.
In order to avoid the singularity condition, the working state of
the system must be restricted \cite{alasty2007nonlinear}. In order
to improve the robustness property of classical feedback linearization,
robust feedback linearization \cite{guillard2000robust,de2018robust,huang2016input}
has been developed. It brings the original nonlinear system into the
Jacobian approximation around the origin rather than in the Brunowski
form. In this case, the original nonlinear system is only partially
transformed, which allows to preserve the good robustness property
obtained by a linear control law which it is associated with. However,
other disadvantages of feedback linearization still remain in robust
feedback linearization. There are also some other linearization methods,
for example, the feedback linearization using Gaussian processes to
predict the unknown functions \cite{umlauft2017feedback}.

To complement the existing linearization methods, a new linearization
method named \emph{state compensation linearization}, which is based
on the additive state decomposition, is proposed in this paper. The
idea of compensation has already existed in active disturbance rejection
control (ADRC) \cite{han2009pid},\cite{chen2015disturbance}, which
can also be used to accomplish linearization, although it mainly aims
to compensate for disturbances. The nonlinear part of a dynamical
system has been forced to be a part of lumped disturbances, which
can be compensated via a feedforward strategy based on the estimates
from an extended state observer (ESO) \cite{Qi2022Problems},\cite{8424483}.
The final linearized system is a cascade integral plant. The significant
feature of ADRC is that it requires minimum information about a dynamic
system except for the relative degree of the system. However, for
most physical systems, the nonlinear dynamics may be known or at least
partially known. The control effect can be improved if the known nonlinear
dynamics could be exploited in design. Besides, how to select the
relative degree for non-minimum phase (NMP) systems is still an open
problem in ADRC \cite{chen2015disturbance}.

\emph{Additive state decomposition }\cite{quan2015additive},\cite{wei2016output}
is\emph{ }a kind of generalization of the superposition principle
in linear systems to nonlinear systems. As a task-oriented problem
decomposition, additive state decomposition\emph{ }can reduce the
complexity of a task, which is the addition of several\emph{\ }subtasks.
The essential idea of additive state decomposition can be expressed
by the following equation{ 
\begin{equation}
{A=B+}\underset{C}{\underbrace{\left(A-B\right)}}
\end{equation}
} where $A$ is the original system/task/problem, $B$ is the primary
system/task/problem, and $C=A-B$ represents the secondary system/task/problem.
It is obvious that the equality holds. In order to accomplish linearization,
the linear system (primary system ${B}$) can be obtained by compensating
for the nonlinearity of the original system using the secondary system
$C$. This can be described as 
\begin{equation}
{B=A-}C.
\end{equation}
The secondary system $C$ can be designed as a system without any
uncertainties. On the other hand, because the sum of the primary system
and the secondary system is the original system, all the uncertainties,
such as the higher-order unmodeled dynamics of the real system, are
left in the primary system (the linear system). Based on the proposed
state compensation linearization method, a \textit{state-compensation-linearization-based
control framework} is then established. The basic idea is to decompose
an original system into several subsystems taking charge of simpler
subtasks. Then one can design controllers for these subtasks and finally
combine them to achieve the original control task. Such control methods
were also proposed in our previous work in \cite{quan2015additive},\cite{quan2014output}.
However, the method was only applied to the control problem for some
special systems. In this paper, the state-compensation-linearization-based
control framework is proposed for more general systems.

The main contributions of this paper are twofold.

(i) A new linearization method called state compensation linearization
is proposed. It is a linearization method complementary to the Jacobian
linearization method and the feedback linearization method. Unlike
the existing two linearization methods, the state compensation linearization
method uses a \emph{compensation strategy}. The linear system (primary
system) can be obtained by compensating for the nonlinearity of the
original system using the secondary system and its controller. This
makes state compensation linearization flexible and combines the advantages
of Jacobian linearization and feedback linearization.

(ii) A state-compensation-linearization-based control framework is
established. Because the primary system is linear, the classical frequency-domain
and time-domain methods for linear systems can be applied. In this
new framework, the nonlinearity of the original system is \textit{compensated}
by the nonlinear secondary system and its stabilizing controller.
Then, the original stabilization/tracking problem for the nonlinear
system is simplified to a stabilization/tracking problem for the linear
primary system, and different control methods are allowed to be applied
to the system.

The remaining parts of this paper are organized as follows. In Section
\ref{Preliminaries}, necessary preliminaries are introduced. Next,
Section \ref{CLMethod} details a new linearization method called
state compensation linearization, which is followed by state compensation
linearization applications to control design in Section \ref{CL based C}.
Three illustrative examples are provided in Section \ref{sec:Simulation-study}
to show the effectiveness of the proposed method. Section \ref{Conclusions}
concludes the paper.

\section{Preliminaries: Additive state decomposition}

\label{Preliminaries}

\emph{Additive state decomposition} has already appeared in some other
domains. The additive state decomposition mentioned here is limited
to the meaning in the author's paper \cite{quan2015additive}. Consider
the following class of differential dynamic systems 
\begin{equation}
A:\left\{ \begin{array}{l}
\mathbf{\dot{x}}=\mathbf{f}\left({\mathbf{x,}}\mathbf{u}\right)+\mathbf{d},\mathbf{x}\left(0\right)=\mathbf{x}_{0}\\
\mathbf{y}=\mathbf{h}\left({\mathbf{x}}\right)
\end{array}\right.\label{Dif_Orig_Sys}
\end{equation}
where $\mathbf{x}\in\mathbb{R}^{n}$ is the state, $\mathbf{u}\in\mathbb{R}^{m}$
is the input, $\mathbf{y}\in\mathbb{R}^{p}$ is the output, $\mathbf{d}\in\mathbb{R}^{n}$
denotes the disturbance, functions $\mathbf{f:\mathbb{R}}^{n}\mathbf{\times\mathbb{R}}^{m}\rightarrow\mathbb{R}^{n}$
and $\mathbf{\mathbf{h}:\mathbb{R}}^{n}\rightarrow\mathbb{R}^{m}$.
There is no restriction on functions $\mathbf{f}$ and $\mathbf{h}$.
For the system (\ref{Dif_Orig_Sys}), we make

\textbf{Assumption 1}. For a given external input $\mathbf{u}$, the
system (\ref{Dif_Orig_Sys}) with the initial value $\mathbf{x}_{0}$
has a unique solution $\mathbf{x}^{*}$ on $\left[0,\infty\right)$.

Noteworthy, most systems in practice satisfy\emph{ Assumption 1},\emph{
}i.e., the uniqueness of solutions\emph{.}

First, a \emph{primary system} is brought in, having the same dimension
as the original system: 
\begin{equation}
B:\left\{ \begin{array}{l}
\mathbf{\dot{x}}_{\text{p}}=\mathbf{f}_{\text{p}}\left({\mathbf{x}_{\text{p}},}\mathbf{u}_{\text{p}}\right)+\mathbf{d},\mathbf{x}_{\text{p}}\left(0\right)=\mathbf{x}_{\text{p},0}\\
\mathbf{y}_{\text{p}}=\mathbf{h}_{\text{p}}\left({\mathbf{x}_{\text{p}}}\right)
\end{array}\right.\label{Dif_Pri_Sys}
\end{equation}
where {$\mathbf{x}_{\text{p}}$}$\in
\mathbb{R}
^{n},\mathbf{u}_{\text{p}}\in
\mathbb{R}
^{m}$, $\mathbf{y}_{\text{p}}\in
\mathbb{R}
^{p}$, $\mathbf{\mathbf{f}_{\text{p}}:\mathbb{R}}^{n}\mathbf{\times\mathbb{R}}^{m}\rightarrow\mathbb{R}^{n}$
and $\mathbf{\mathbf{h}_{\text{p}}:\mathbb{R}}^{n}\rightarrow\mathbb{R}^{m}$.
Functions $\mathbf{f}_{\text{p}}$ and $\mathbf{h}_{\text{p}}$ are
determined by the designer. From the original system (\ref{Dif_Orig_Sys})
and the primary system (\ref{Dif_Pri_Sys}), the following \emph{secondary
system} is derived: 
\[
C=A-B:\left\{ \begin{array}{l}
\mathbf{\dot{x}}-\mathbf{\dot{x}}_{\text{p}}=\mathbf{f}\left({\mathbf{x,}}\mathbf{u}\right)-\mathbf{f}_{\text{p}}\left({\mathbf{x}_{\text{p}},}\mathbf{u}_{\text{p}}\right)\\
\mathbf{y}-\mathbf{y}_{\text{p}}=\mathbf{h}\left({\mathbf{x}}\right)-\mathbf{h}_{\text{p}}\left({\mathbf{x}_{\text{p}}}\right).
\end{array}\right.
\]
New variables {$\mathbf{x}_{\text{s}}$}$\in
\mathbb{R}
^{n},\mathbf{u}_{\text{s}}\in
\mathbb{R}
^{m},$\textbf{\ }$\mathbf{y}{_{\text{s}}}\in
\mathbb{R}
^{p}$ are defined as follows 
\begin{equation}
{\mathbf{x}_{\text{s}}}\triangleq{\mathbf{x}}-{\mathbf{x}_{\text{p}},}\mathbf{u}{_{\text{s}}}\triangleq\mathbf{u}-\mathbf{u}{_{\text{p}},}\mathbf{y}{_{\text{s}}}\triangleq\mathbf{y}-\mathbf{y}_{\text{p}}.\label{Gen_RelationPS}
\end{equation}
Then, the secondary system can be further written as follows 
\begin{align}
\mathbf{\dot{x}}_{\text{s}} & =\mathbf{f}\left({\mathbf{x}_{\text{p}}}+{\mathbf{x}_{\text{s}},}\mathbf{u}_{\text{p}}+\mathbf{u}{_{\text{s}}}\right)-\mathbf{f}_{\text{p}}\left({\mathbf{x}_{\text{p}},\mathbf{u}_{\text{p}}}\right),\mathbf{x}_{\text{s}}\left(0\right)=\mathbf{x}_{\text{s,0}}\nonumber \\
\mathbf{y}_{\text{s}} & =\mathbf{h}\left({\mathbf{x}_{\text{p}}}+{\mathbf{x}_{\text{s}}}\right)-\mathbf{h}_{\text{p}}\left({\mathbf{x}_{\text{p}}}\right).\label{Dif_Sec_Sys}
\end{align}
From the definition (\ref{Gen_RelationPS}), it follows 
\begin{equation}
{\mathbf{x}}={\mathbf{x}_{\text{p}}}+{\mathbf{x}_{\text{s}},}\mathbf{u}=\mathbf{u}{_{\text{p}}}+\mathbf{u}{_{\text{s}},\mathbf{y}}={\mathbf{y}_{\text{p}}}+{\mathbf{y}_{\text{s}}.}\label{Relation}
\end{equation}

Now we can state

\textbf{Lemma 1.} Under \emph{Assumption 1}, suppose $\mathbf{x}_{\text{p}}$
and $\mathbf{x}_{\text{s}}$ are the solutions of the system (\ref{Dif_Pri_Sys})
and system (\ref{Dif_Sec_Sys}) respectively, and the initial conditions
of (\ref{Dif_Orig_Sys}), (\ref{Dif_Pri_Sys}) and (\ref{Dif_Sec_Sys})
satisfy $\mathbf{x}_{0}=\mathbf{x}_{\text{p},0}+\mathbf{x}_{\text{s,0}}$.
Then ${\mathbf{x}}={\mathbf{x}_{\text{p}}}+{\mathbf{x}_{\text{s}}}$.

\textit{\emph{Proof.}}\textit{ }\textit{\emph{See reference \cite{quan2009additive}.}}

Additive state decomposition is applicable not only to the above class
of differential dynamic systems but also to almost any dynamic or
static systems. As for a concrete decomposition process, it requires
a specific analysis for a specific problem.

\section{A new linearization method: State compensation linearization}

\label{CLMethod} In this section, two existing linearization methods
are recalled firstly. Then, as a new linearization method, state compensation
linearization is proposed. Finally, the comparison of the three linearization
methods is summarized.

For clarification, the following general nonlinear system is considered
\begin{equation}
\mathbf{\dot{x}=f}\left(\mathbf{x,u}\right)+\mathbf{d}\mathbf{,}\text{ }\mathbf{x}\left(0\right)=\mathbf{x}_{0}\label{generalnonlinear}
\end{equation}
where $\mathbf{x}\in\mathbb{R}^{n}$ is the state, $\mathbf{u}\in\mathbb{R}^{m}$
is the input, $\mathbf{d}\in\mathbb{R}^{n}$ denotes the disturbance,
$\mathbf{x}_{0}$ is the initial state, and $\mathbf{f:\mathbb{R}}^{n}\mathbf{\times\mathbb{R}}^{m}\rightarrow\mathbb{R}^{n}$.
For simplicity, let $\mathbf{f}\left(\mathbf{0,0}\right)=\mathbf{0,}$
which implies that the origin is the equilibrium point of the nominal
system (\ref{generalnonlinear}) with $\mathbf{d}=\mathbf{0}$.

\textbf{Assumption 2}. The state $\mathbf{x}$ is measurable.

\textbf{Assumption 3}. The function $\mathbf{f}\left(\mathbf{x,u}\right)$
is continuously differentiable relative to $\mathbf{x}$ and $\mathbf{u}$
around the equilibrium point.

\subsection{Existing linearization methods}

There are currently two main linearization methods: \emph{Jacobian
linearization (i.e., linearization by Taylor's expansion)} and \emph{feedback
linearization}.

\subsubsection{Jacobian linearization}

Jacobian linearization is based on the equilibrium point of the considered
system. Considering that, when $\mathbf{d}=\mathbf{0}$, the equilibrium
point of system (\ref{generalnonlinear}) is the origin, the nonlinear
system is linearized as 
\begin{equation}
\mathbf{\dot{x}=A}_{1}\mathbf{x}+\mathbf{B}_{1}\mathbf{u}\label{Taylorlinear}
\end{equation}
where $\mathbf{A}_{1}=\left.\frac{\partial\mathbf{f}}{\partial\mathbf{x}}\right\vert _{\mathbf{x}=\mathbf{0,u}=\mathbf{0}}\in\mathbf{
\mathbb{R}
}^{n\times n},\mathbf{B}_{1}=\left.\frac{\partial\mathbf{f}}{\partial\mathbf{u}}\right\vert _{\mathbf{x}=\mathbf{0,u}=\mathbf{0}}\in\mathbf{
\mathbb{R}
}^{n\times m}$. The detailed conditions for Jacobian linearization can be seen in
\cite{slotine1991applied}. Then, a controller is often designed as
$\mathbf{u=L}_{1}\left(\mathbf{x}\right)$, where $\mathbf{L}_{1}\left(\mathbf{\cdot}\right)$
is a linear function, and the feedback state $\mathbf{x}$ is that
of the real system (\ref{generalnonlinear}).

Jacobian linearization is the most widely-used method to deal with
nonlinear systems in practice. It is simple and intuitive. Moreover,
this linearization method is easy-to-use because it can be applied
to most nonlinear functions and the linearization process can be implemented
by a program automatically. After Jacobian linearization, the physical
meaning of states remains the same as before so that it can be easily
understood. However, this method has only a local property. This implies
that linearization may cause worse performance, or even instability
of the resulting closed-loop system when the states are far away from
the equilibrium points (trim states). In order to reduce approximation
error, the number of equilibrium points required for the entire space
has to be increased. It may be particularly large. As a result, this
will result in a large number of linear models \cite{nichols1993gain}.
Then, the controller gains need to be scheduled according to the linear
models \cite{leith2000survey}. Trajectory linearization is a special
kind of Jacobian linearization that linearizes the system at every
point on the nominal trajectory \cite{liu2008omni}. In essence, the
obtained linearized systems are time-varying, and time-varying analysis
and control design are necessary.

\subsubsection{Feedback linearization}

First, it should be pointed out that the feedback linearization method
is only applicable to precisely known systems without disturbances
and uncertainties. Thus, it is assumed that $\mathbf{d}=\mathbf{0}$
in the following discussion. Note that the following affine nonlinear
systems are often considered in feedback linearization 
\begin{equation}
\mathbf{\dot{x}=}\mathbf{f}\left(\mathbf{x}\right)+\mathbf{g}\left(\mathbf{x}\right)\mathbf{u}\label{AfineN}
\end{equation}
where $\mathbf{f}:\mathbf{\mathbb{R}}^{n}\rightarrow\mathbb{R}^{n}$
and $\mathbf{g}:\mathbf{\mathbb{R}}^{n}\rightarrow\mathbb{R}^{n\times m}$
are smooth vector fields. The detailed conditions for feedback linearization
can been seen in \cite{slotine1991applied}. Find a state transformation
$\mathbf{z=T}\left(\mathbf{x}\right)\in
\mathbb{R}
^{n}$ and an input transformation $\mathbf{u=u}\left(\mathbf{z,v}\right)\in
\mathbb{R}
^{m}$ so that the system (\ref{generalnonlinear}) can be converted into
the following form \cite{slotine1991applied} 
\begin{equation}
\mathbf{\dot{z}=A}_{2}\mathbf{z+B}_{2}\mathbf{v}\label{FeedbackLinearization}
\end{equation}
where $\mathbf{A}_{2}\in
\mathbb{R}
^{n\times n},$ $\mathbf{B}_{2}\in
\mathbb{R}
^{n\times m}$. For system (\ref{FeedbackLinearization}), the controller is designed
as $\mathbf{v=L}_{2}\left(\mathbf{z}\right)$. The feedback control
is based on the variable $\mathbf{z}$ after the coordinate transformation.
Finally, the actual control input can be obtained by $\mathbf{u=u}\left(\mathbf{z,v}\right)$.

After feedback linearization, system (\ref{FeedbackLinearization})
no longer has multiple equilibrium points and has a linear property
in the entire transformed space. However, the coordinate transformation
is complicated. There are often many situations where feedback linearization
is difficult to apply: (i) the exact model of the nonlinear system
is not available, and no robustness is guaranteed in the presence
of parameter uncertainties or unmodeled dynamics; (ii) the control
law is not defined at singularity points; (iii) it is not easy to
use for NMP systems, which may lead to unstable control; (iv) the
physical meaning of new states differs from that of the original system,
which makes the controller parameter adjustment difficult.

To avoid the disadvantages of the two existing linearization methods,
a new linearization method named \emph{state compensation linearization},
which is based on the additive state decomposition, is proposed as
a kind of complementation in the next subsection.

\subsection{State compensation linearization as a kind of complementation}

\label{subsec:State-compensation-linearization}

In this part, the additive state decomposition is applied to linearization,
which leads to a new linearization method named state compensation
linearization. For the original nonlinear system (\ref{generalnonlinear}),
denoted as $A$, the primary system is chosen as 
\begin{equation}
\mathbf{\dot{x}}_{\text{p}}=\mathbf{A}_{1}\mathbf{x}_{\text{p}}+\mathbf{B}_{1}\mathbf{u}_{\text{p}}+\mathbf{d},\mathbf{x}_{\text{p}}\left(0\right)=\mathbf{x}_{0}\label{CL_PSys}
\end{equation}
which is a linear system and denoted as $B$. Here, $\mathbf{A}_{1}=\left.\frac{\partial\mathbf{f}}{\partial\mathbf{x}}\right\vert _{\mathbf{x}=\mathbf{0,u}=\mathbf{0}}\in\mathbf{
\mathbb{R}
}^{n\times n},\mathbf{B}_{1}=\left.\frac{\partial\mathbf{f}}{\partial\mathbf{u}}\right\vert _{\mathbf{x}=\mathbf{0,u}=\mathbf{0}}\in\mathbf{
\mathbb{R}
}^{n\times m}$, which are the same as in (\ref{Taylorlinear}). Then, the secondary
system is determined by the original system and the primary system
using the rule (\ref{Dif_Sec_Sys}), namely $C=A-B$, that 
\begin{equation}
\mathbf{\dot{x}}_{\text{s}}=\mathbf{f}\left(\mathbf{x},\mathbf{u}\right)-\mathbf{A}_{1}\mathbf{x}_{\text{p}}-\mathbf{B}_{1}\mathbf{u}_{\text{p}},\mathbf{x}_{\text{s}}\left(0\right)=\mathbf{0},\label{CL_SSys}
\end{equation}
which can be then written as 
\begin{equation}
\mathbf{\dot{x}}_{\text{s}}=\mathbf{f}\left(\mathbf{x},\mathbf{u}\right)+\mathbf{A}_{1}\left(\mathbf{x}_{\text{s}}-\mathbf{x}\right)+\mathbf{B}_{1}\left(\mathbf{u}_{\text{s}}-\mathbf{u}\right),\mathbf{x}_{\text{s}}\left(0\right)=\mathbf{0}\label{CL_SSys-1}
\end{equation}
where 
\begin{equation}
{\mathbf{x}}={\mathbf{x}_{\text{p}}}+{\mathbf{x}_{\text{s}},}\mathbf{u}=\mathbf{u}{_{\text{p}}}+\mathbf{u}{_{\text{s}}.}\label{eq:relation}
\end{equation}

\begin{figure}[tbph]
\begin{centering}
\includegraphics[scale=0.7]{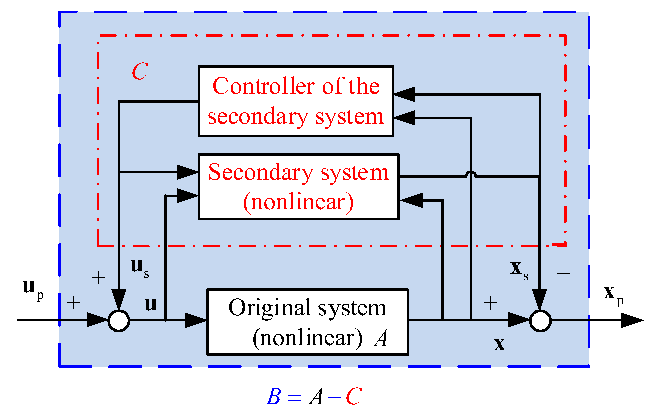} 
\par\end{centering}
\caption{The construction diagram of the linear primary system from the input
$\mathbf{u}_{\text{p}}$ to the state $\mathbf{x}_{\text{p}}$. The
primary system can be obtained by subtracting the complementary secondary
system from the real system.}
\label{CL} 
\end{figure}

In order to accomplish linearization, the linear system (\ref{CL_PSys})
($B=A-C$) can be obtained by compensating for the nonlinearity of
the original system ($A$) using the secondary system ($C$) and its
controller $\mathbf{u}_{\text{s}}=\mathbf{L}_{4}\left(\mathbf{x},\mathbf{x}_{\text{s}}\right)$
as shown in Fig. \ref{CL}. The system from the input $\mathbf{u}_{\text{p}}$
to the state $\mathbf{x}_{\text{p}}$ is linear. In the subsystems
(\ref{CL_PSys}) and (\ref{CL_SSys}), the secondary system (\ref{CL_SSys})
is a virtual exact system, which exists only in design. On the other
hand, since the sum of the primary system and the secondary system
is the original system, all the uncertainties, such as the higher-order
unmodeled dynamics of the real system, are left in the primary system.
After state compensation linearization, the controller for the linear
primary system (\ref{CL_PSys}) can be designed as $\mathbf{u}_{\text{p}}=\mathbf{L}_{3}\left(\mathbf{x}_{\text{p}}\right)$.

Based on the above content, the following theorem about the conditions
for state compensation linearization can be stated.

\textbf{Theorem 1}.\textit{ }\textit{\emph{System }}(\ref{generalnonlinear})\textit{\emph{
can be linearized by }}state compensation linearization\textit{\emph{
if}}

\textit{\emph{i). }}in the entire space $\mathbf{x}\in\mathbb{R}^{n},\mathbf{u}\in\mathbb{R}^{m}$,
\emph{Assumption 1 }is satisfied.

ii). \textit{\emph{it can be linearizable by }}Jacobian linearization\textit{\emph{.}}

\textit{\emph{Proof. First, if Condition i) is met, according to }}\textit{Lemma
1}\textit{\emph{, }}additive state decomposition\textit{\emph{ holds.
Next, if system }}(\ref{generalnonlinear}) can be\emph{ }\textit{\emph{linearizable
by}}\textit{ }Jacobian linearization, then the matrices $\mathbf{A}_{1},\mathbf{B}_{1}$
can be obtained, and thus the primary system (\ref{CL_PSys}) can
be determined. Then, \textit{\emph{according to}} additive state decomposition,
the secondary system (\ref{CL_SSys}) is determined by the original
system and the primary system using the rule (\ref{Dif_Sec_Sys}),
namely $C=A-B$. Next, one can design a controller $\mathbf{u}_{\text{s}}=\mathbf{L}_{4}\left(\mathbf{x},\mathbf{x}_{\text{s}}\right)$
to stabilize the secondary system. Consequently, by compensating for
the nonlinearity of the original system using the secondary system
and its controller, the linear system (i.e., the primary system (\ref{CL_PSys}))
can be obtained. \emph{$\square$}

State compensation linearization no longer has multiple equilibrium
points. Moreover, in state compensation linearization, ${\mathbf{x}}={\mathbf{x}_{\text{p}}}+{\mathbf{x}_{\text{s}},}$
there are no coordinate transformation or unit transformation. In
some sense, $\mathbf{x}_{\text{p}}$ can be considered as a part of
$\mathbf{x}$, and they have the same dimension and unit. Thus, the
physical meaning of the states remains the same after state compensation
linearization. This property of state compensation linearization is
similar to that of Jacobian linearization. In contrast, the new states
obtained from feedback linearization have different physical meanings
with the original states because of coordinate transformation.

\textbf{Remark 1}. In state compensation linearization, one keeps
both `linear primary' and `nonlinear secondary' component systems.
Furthermore, it is not necessarily assumed that the secondary system
is `small' in some sense. The secondary system may be significant
when the original system has a strong nonlinearity. However, what
we often do is to stabilize the secondary system so that it does become
`small' eventually, leaving one free to work with the primary system.
$A$ and $B$ are said to be equivalent if the nonlinear secondary
system $A-B$ can be stabilized by a controller. Then, checking the
equivalence means checking whether the system $A-B$ can be stabilized.
In that case a nonlinear system $A$ is equivalent to the linear system
$B$ for which $A-B$ can be stabilized. The difference between state
compensation linearization and Jacobian linearization is that state
compensation linearization does not throw the `error' system in Jacobian
linearization away.

\textbf{Remark 2}. In state compensation linearization, one often
needs to construct a stabilizing controller for a nonlinear secondary
system. Although designing a stabilizing controller may cost some
effort, the benefit is also attractive. The nonlinearity in the secondary
system can be better considered, and thus better performance can be
expected. Besides, because all the uncertainties and the initial value
are considered in the primary system, the stabilizing control problem
for the secondary system is a stabilizing control problem of an exact
nonlinear system with zero initial states, which is often not complicated.

\subsection{Comparison of the three linearization methods}

Although the linearized systems (\ref{CL_PSys}) and (\ref{Taylorlinear})
have the same $\mathbf{A}_{1},\mathbf{B}_{1}$, the state and input
in (\ref{CL_PSys}) are $\mathbf{x}_{\text{p}}$ and $\mathbf{u}_{\text{p}}$
instead of the real state and input $\mathbf{x}$ and $\mathbf{u}$
in (\ref{Taylorlinear}). 

Under some specific conditions, state compensation linearization can
be transformed into Jacobian linearization. Considering the original
system (\ref{generalnonlinear}) with $\mathbf{d}=\mathbf{0}$. In
state compensation linearization, let $\ensuremath{\mathbf{x}_{\text{s}}=\mathbf{0}}$,
$\mathbf{u}_{\text{s}}=\mathbf{0}$, then the resulting linearized
system is the same as that of Jacobian linearization, namely $\mathbf{\dot{x}=A}_{1}\mathbf{x}+\mathbf{B}_{1}\mathbf{u}$.

State compensation linearization seems similar to Jacobian linearization.
The difference is that the latter drops the higher-order small nonlinear
terms directly and Taylor's expansion holds locally, which determines
its approximate and local properties. In contrast, state compensation
linearization uses a \emph{compensation strategy}, which makes state
compensation linearization extend to the whole space and become a
kind of accurate linearization method. This property is similar to
that of feedback linearization. Thus, state compensation linearization
combines the advantages of the two linearization methods in some sense.

Along with state compensation linearization, an original control task
can be assigned to the linearized primary system and nonlinear secondary
system by dividing the original control task into two simpler control
tasks. From this perspective, additive state decomposition\emph{ }is
also a \textit{task-oriented problem decomposition} that offers a
general way to decompose an original system into several subsystems
in charge of simpler subtasks. Thanks to task decomposition, one can
design controllers for simpler subtasks more easily than designing
a controller for the original task. These controllers can finally
be combined to accomplish the original control task. In the following,
state compensation linearization applications to nonlinear control
design are discussed.

\section{State compensation linearization applications to nonlinear control
design}

\label{CL based C}

Based on the state compensation linearization method in the previous
section, a state-compensation-linearization-based control framework
for a class of nonlinear systems is proposed in this subsection. State-compensation-linearization-based
stabilizing control is discussed in the following for illustration.

\subsection{Problem formulation}

Consider the nonlinear system (\ref{generalnonlinear}), state compensation
linearization has been done in Subsection \ref{subsec:State-compensation-linearization}.
Then, the control objective is to make $\mathbf{x}\left(t\right)\rightarrow\mathbf{0}$
or $\mathbf{x}\left(t\right)\rightarrow\mathcal{B}\left(\mathbf{0}_{n\times1},\sigma\right)$\footnote{$\mathcal{B}\left(\mathbf{o},\sigma\right)\triangleq\left\{ \mathbf{x}\in\mathbb{R}^{n}\left\vert \left\Vert \mathbf{x}-\mathbf{o}\right\Vert \leq\sigma\right.\right\} ,$
and $\boldsymbol{\xi}\left(t\right)\rightarrow\mathcal{B}\left(\mathbf{o},\sigma\right)$
signifies $\underset{\mathbf{y}\in\mathcal{B}\left(\mathbf{o},\sigma\right)}{\min}$
$\left\Vert \boldsymbol{\xi}\left(t\right)-\mathbf{y}\right\Vert \rightarrow0.$} as $t\rightarrow\infty$.

\subsection{Observer design}

\label{subsec:Observer-design}

It should be pointed out that the primary system and the secondary
system are virtual systems, which do not exist in practice and just
exist in the analysis. Thus, an observer is necessary to provide the
variables of the primary system and the secondary system for the later
control design. By using the state compensation linearization method,
the unknown disturbance and the initial value are assigned in the
primary system, and the secondary system is an exact nonlinear system.
Thus, the following observer is built based on the secondary system
to observe states \cite{quan2015additive},\cite{quan2014output}.

\textbf{Theorem 2}.\textit{ }\textit{\emph{Suppose that }}$\mathbf{A}_{1}$\textit{\emph{
is stable, and an observer is designed to estimate }}$\mathbf{x}_{\text{p}}$\emph{,
}\textit{\emph{$\mathbf{x}_{\text{s}}$ in (\ref{CL_PSys}) and (\ref{CL_SSys})
as}} 
\begin{align}
\mathbf{\dot{\hat{x}}}_{\text{s}} & =\mathbf{f}\left(\mathbf{x},\mathbf{u}\right)+\mathbf{A}_{1}\left(\mathbf{\hat{x}}_{\text{s}}-\mathbf{x}\right)+\mathbf{B}_{1}\left(\mathbf{u}_{\text{s}}-\mathbf{u}\right),\mathbf{\hat{x}}_{\text{s}}\left(0\right)=\mathbf{0}\label{eq:Obs1}\\
\mathbf{\hat{x}}_{\text{p}} & =\mathbf{x-\hat{x}}_{\text{s}}\label{eq:Obs2}
\end{align}
\textit{\emph{where $\mathbf{\hat{x}}_{\text{s}}$ and}}\emph{ }$\mathbf{\hat{x}}_{\text{p}}$\textit{\emph{
are the estimation of $\mathbf{x}_{\text{s}}$ and}}\emph{ }$\mathbf{x}_{\text{p}}$\textit{\emph{.
Then }}$\hat{\mathbf{x}}_{\text{p}}=\mathbf{x}_{\text{p}}$,\textit{
}$\mathbf{\hat{x}}_{\text{s}}=\mathbf{x}_{\text{s}}$\textit{\emph{.}}

\textit{Proof. }Subtracting (\ref{eq:Obs1}) from (\ref{CL_SSys})
results in $\dot{\tilde{\mathbf{x}}}_{\text{s}}=\mathbf{A}_{1}\tilde{\mathbf{x}}_{\text{s}}$
and $\tilde{\mathbf{x}}_{\text{s}}\left(0\right)=\mathbf{0}$, where
$\tilde{\mathbf{x}}_{\text{s}}\triangleq\mathbf{x}_{\text{s}}-\hat{\mathbf{x}}_{\text{s}}$.
Then, considering that $\mathbf{A}_{1}$ is stable, $\tilde{\mathbf{x}}_{\text{s}}=\mathbf{0}$,
which implies that $\mathbf{\hat{x}}_{\text{s}}=\mathbf{x}_{\text{s}}$.
Consequently, by (\ref{eq:relation}), it can be obtained that $\mathbf{\hat{x}}_{\text{p}}=\mathbf{x-\hat{x}}_{\text{s}}=\mathbf{x}_{\text{p}}$.
$\square$

\textbf{Remark 3. }The designed observer is an open-loop observer,
in which $\mathbf{A}_{1}$ must be stable. Otherwise, state feedback
is needed to obtain a stable $\mathbf{A}_{1}$. If\textbf{ }a closed-loop
observer is adopted here, it will be rather difficult to analyze the
stability of the closed-loop system consisting of a nonlinear controller
and the observer because separation principle does not hold in nonlinear
systems.

\textbf{Remark 4. }The initial values $\mathbf{x}_{\text{s}}\left(0\right)$
and $\mathbf{\hat{x}}_{\text{s}}\left(0\right)$ are both assigned
as $\mathbf{x}_{\text{s}}\left(0\right)=\mathbf{0}$ and $\mathbf{\hat{x}}_{\text{s}}\left(0\right)=\mathbf{0}$.
If there is an initial value measurement error, it will be assigned
to and considered in the primary system in the form of disturbances.

\subsection{State-compensation-linearization-based control framework}

Based on state compensation linearization, the original system (\ref{generalnonlinear})
is decomposed into two systems: an LTI system including all disturbances
as the primary system, together with the secondary system whose equilibrium
point is the origin. Since the states of the primary system and the
secondary system can be observed, the original stabilizing problem
for the system (\ref{generalnonlinear}) is correspondingly decomposed
into two problems: a stabilizing problem for the nonlinear secondary
system (for the nonlinearity compensation objective) and a stabilizing
problem for an LTI primary system (for the control objective). Because
the primary system is linear, the classical frequency-domain and time-domain
methods can be applied, such as lead-lag compensation. Therefore,
the original problem turns into two simpler subproblems. The final
objective is achieved once the objectives of the two subproblems are
achieved. 
\begin{itemize}
\item \textbf{Problem 1} Consider the secondary system (\ref{CL_SSys}).
Design the secondary controller as 
\begin{equation}
\mathbf{u}_{\text{s}}=\mathbf{L}(\mathbf{x},\mathbf{x_{\text{s}}})\label{eq:Controller_Prob2-1}
\end{equation}
such that $\mathbf{x}_{\text{s}}\left(t\right)\rightarrow\mathbf{0}$
as $t\rightarrow\infty$, where $\mathbf{L}(\cdot,\cdot)$ is a nonlinear
function. 
\end{itemize}
After the compensation of the secondary system and its controller
(\ref{eq:Controller_Prob2-1}), one can get a linear system (primary
system). Next, just a stabilizing problem for a linear system needs
to be considered. 
\begin{itemize}
\item \textbf{Problem 2} Consider the primary system (\ref{CL_PSys}). Design
the primary controller as 
\begin{equation}
\mathbf{u}_{\text{p}}=\mathcal{L}^{-1}(\mathbf{H}(s)\mathbf{K}\mathbf{x_{\text{p}}}(s))\label{eq:Controller_Prob1-1}
\end{equation}
such that $\mathbf{x}_{\text{p}}\left(t\right)\rightarrow\mathbf{0}$
or $\mathbf{x}_{\text{p}}\left(t\right)\rightarrow\mathcal{B}\left(\mathbf{0}_{n\times1},\sigma\right)$
as $t\rightarrow\infty$, where $\text{\textbf{H}(s)}$ is a stable
transfer function, $\mathcal{L}^{-1}$ denotes the inverse Laplace
transformation, and state feedback matrix $\mathbf{K\in\mathbb{R}^{\text{\ensuremath{m\times n}}}}$. 
\end{itemize}
\textbf{Remark 5}. In the controller (\ref{eq:Controller_Prob2-1}),
the state feedback matrix $\mathbf{K}$ can be determined by the linear
quadratic regulator (LQR) method or the ``PID Tuner App'' in the
Matlab software. Besides, one can use some kind of robust methods
to construct a controller for the primary system.

The controller of the original system can be obtained by combing the
controllers of the primary system and the secondary system. With the
solutions to the two problems in hand, we can state

\textbf{Theorem 3.} For system (\ref{generalnonlinear}) under \emph{Assumptions
1-3}, suppose (i)\textit{\emph{ }}$\mathbf{A}_{1}$\textit{ is stable;}\emph{
}(ii)\emph{ Problems 1-2} are solved; (iii) the controller for system
(\ref{generalnonlinear}) is designed as 
\begin{align}
\mathbf{\dot{\hat{x}}}_{\text{s}} & =\mathbf{f}\left(\mathbf{x},\mathbf{u}\right)+\mathbf{A}_{1}\left(\mathbf{\hat{x}}_{\text{s}}-\mathbf{x}\right)+\mathbf{B}_{1}\left(\mathbf{u}_{\text{s}}-\mathbf{u}\right)\nonumber \\
\hat{\mathbf{x}}_{\text{p}} & =\mathbf{x}-\mathbf{\hat{x}}_{\text{s}}\nonumber \\
\mathbf{u} & \mathcal{=L}^{-1}(\mathbf{H}(s)\mathbf{K}\hat{\mathbf{x}}_{\text{p}}(s))+\mathbf{L}(\mathbf{x},\mathbf{\hat{x}}_{\text{s}}).\label{eq:Controller-1}
\end{align}
Then, the state of system (\ref{generalnonlinear}) satisfies $\mathbf{x}\left(t\right)\rightarrow\mathbf{0}$
or $\mathbf{x}\left(t\right)\rightarrow\mathcal{B}\left(\mathbf{0}_{n\times1},\sigma\right)$
as $t\rightarrow\infty$.

Proof. It is easy to obtain that $\mathbf{\hat{x}}_{\text{p}}=\mathbf{x}_{\text{p}}$
and $\mathbf{\hat{x}}_{\text{s}}=\mathbf{x}_{\text{s}}$ with the
observer (\ref{eq:Obs1})-(\ref{eq:Obs2}). Then, since \textit{Problem
1} is solved, $\mathbf{x}_{\text{s}}\left(t\right)\rightarrow\mathbf{0}$
as $t\rightarrow\infty$. Furthermore, $\mathbf{x}_{\text{p}}\left(t\right)\rightarrow\mathbf{0}$
or $\mathbf{x}_{\text{p}}\left(t\right)\rightarrow\mathcal{B}\left(\mathbf{0}_{n\times1},\sigma\right)$
as $t\rightarrow\infty$ with \textit{Problem 2} solved. By additive
state decomposition, we have $\mathbf{x=\mathbf{x}}_{\text{p}}+\mathbf{x}_{\text{s}}\rightarrow\mathbf{0}$
or $\mathbf{x}\left(t\right)\rightarrow\mathcal{B}\left(\mathbf{0}_{n\times1},\sigma\right)$
as $t\rightarrow\infty$. $\square$

\begin{figure}[ptbh]
\begin{centering}
\includegraphics[scale=0.7]{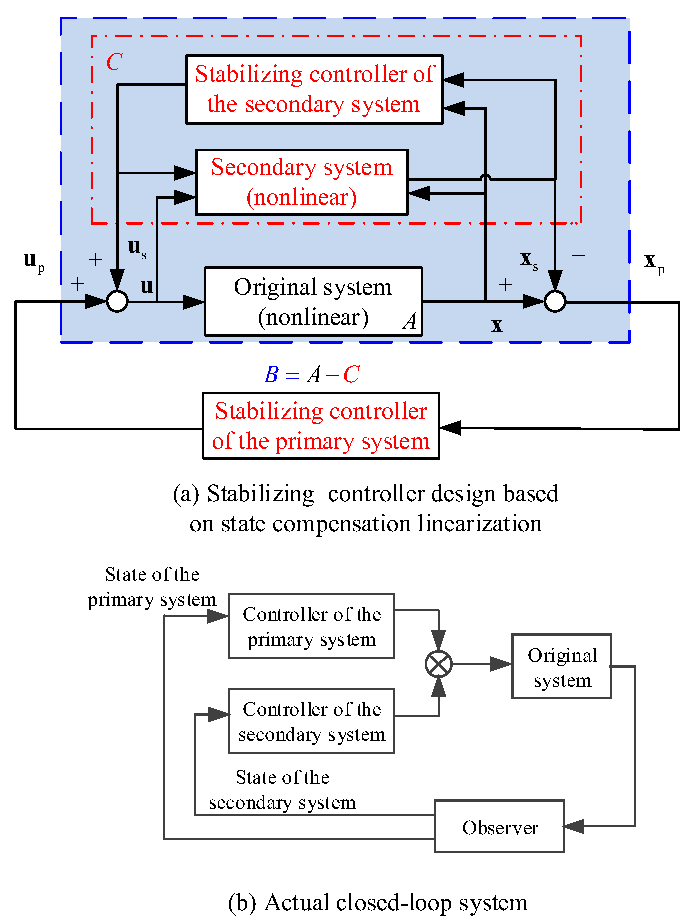} 
\par\end{centering}
\caption{The closed-loop control diagram based on the state-compensation-linearization-based
stabilizing control for a class of nonlinear systems}
\label{closed-loop-s} 
\end{figure}

In this new framework, the original problem is simplified, and different
control methods are allowed to be applied to the systems. The closed-loop
system is shown in Fig. \ref{closed-loop-s}. The controller based
on state compensation linearization includes three parts: the controller
of the primary system, the controller of the secondary system and
the observer. In practice, the primary system and the secondary system
do not need to have physical status. They are only the models existing
in the design. 

\textbf{Remark 6}. The basic idea of the state-compensation-linearization-based
tracking control framework is similar to that of the state-compensation-linearization-based
stabilizing control framework. The most significant difference is
that a tracking control problem rather than a stabilizing control
problem needs to be specified to the linear primary system.

\section{Simulation study}

\label{sec:Simulation-study}

In the following, three examples\textit{ }are given to illustrate
and compare different linearization methods for control design.

\subsection{Example 1: A bilinear system}

Consider a bilinear system as follows 
\begin{align}
\dot{x} & =-4x+xu+d\nonumber \\
y & =x\label{Ex1}
\end{align}
where $x,y,u,d\in
\mathbb{R}
$ represent the state, output, input and disturbance, respectively.
It is assumed that $d\neq0$. The tracking control goal is to make
$y\rightarrow y_{\text{d}}$ as $t\rightarrow\infty$.

\subsubsection{Jacobian-linearization-based-control (JLC) and feedback-linearization-based
control (FLC)}

The two existing linearization methods are attempted to linearize
system (\ref{Ex1}) so that the given tracking control goal can be
achieved by using linear control methods. When the Jacobian linearization
method is applied, the first thing is to obtain equilibrium points
of system (\ref{Ex1}). Unfortunately, the stability property of the
nonlinear system (\ref{Ex1}) is affected by the disturbance $d$,
which is unknown. Additionally, $u$ can be arbitrary values when
$x=0$. This implies that equilibrium points of system (\ref{Ex1})
cannot be obtained. Without equilibrium points, it is impossible to
linearize system (\ref{Ex1}) by Taylor's expansion. When the feedback
linearization method is applied, let $\dot{y}=v$. Then, $v=-4x+xu+d$.
There exists a singularity problem that system (\ref{Ex1}) becomes
uncontrollable \cite{tie2014controllability} when $x=0$. Thus, neither
of the two existing linearization methods is suitable for the bilinear
system.

\subsubsection{State-compensation-linearization-based control (SCLC)}

\begin{figure}[tbph]
\begin{centering}
\includegraphics[scale=0.8]{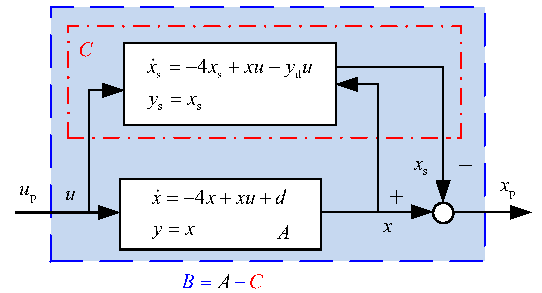} 
\par\end{centering}
\caption{The state compensation linearization structure. The primary system
(\protect\ref{PSys1}) from the input $u_{\text{p}}$ to the state
$x_{\text{p }}$ is linear, which can be obtained by subtracting the
complementary secondary system (\protect\ref{SSys1}) from the original
system (\protect\ref{Ex1}).}
\label{Fig_Ex1-1} 
\end{figure}

\begin{figure}[tbph]
\begin{centering}
\includegraphics[scale=0.65]{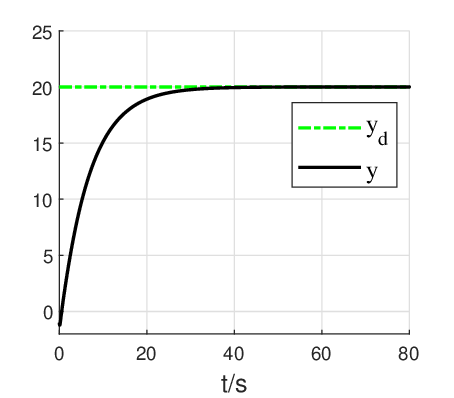} 
\par\end{centering}
\caption{Output response of system (\protect\ref{Ex1}) based on the state
compensation linearization method (Example 1)}
\label{Fig_Ex1-2} 
\end{figure}

Considering the two existing linearization methods are not readily
applicable to system (\ref{Ex1}), the state compensation linearization
method is then adopted. First, a linear primary system ($B$) is designed
as follows 
\begin{align}
\dot{x}_{\text{p}} & =-4x_{\text{p}}+y_{\text{d}}u_{\text{p}}+d,x_{\text{p}}\left(0\right)=x_{0}\nonumber \\
y_{\text{p}} & =x_{\text{p}}\label{PSys1}
\end{align}
where let $u_{\text{p}}=u$, $y_{\text{d}}$ is the desired output.
In the linearization process, because the desired input $u_{\text{d}}$
is unknown, $u_{\text{d}}=0$ is chosen for simplicity. Then, the
following secondary system ($C=A-B$) is determined by subtracting
the primary system (\ref{PSys1}) from the original system (\ref{Ex1})
with rule (\ref{Dif_Sec_Sys}) 
\begin{align}
\dot{x}_{\text{s}} & =-4x_{\text{s}}+xu-y_{\text{d}}u\nonumber \\
 & =\left(u-4\right)x_{\text{s}}+\left(x_{\text{p}}-y_{\text{d}}\right)u\nonumber \\
y_{\text{s}} & =x_{\text{s}}.\label{SSys1}
\end{align}
According to additive state decomposition, it follows 
\begin{equation}
x=x{_{\text{p}}}+x{_{\text{s}},}u=u{_{\text{p}}.}
\end{equation}
Then, the linear system (\ref{PSys1}) ($B=A-C$) can be obtained
by compensating for the nonlinearity of the original system (\ref{Ex1})
($A$) using the secondary system ($C$). The structure is displayed
in Fig. \ref{Fig_Ex1-1}. If the controller for the primary system
(\ref{PSys1}) is well designed, namely $x_{\text{p}}\rightarrow y_{\text{d}}$,
$u\rightarrow u_{\text{d}}$. Then, the secondary system is approaching
the dynamics of $\dot{x}_{\text{s}}=\left(u_{\text{d}}-4\right)x_{\text{s}},$
which implies that $x_{\text{s}}\rightarrow0$ as long as $u_{\text{d}}<4$.
The final linear system (\ref{PSys1}) is independent of system (\ref{SSys1}).
The remaining problem is to design a controller for the linear system
(\ref{PSys1}), for which any linear control methods can be utilized.
Noteworthy, the proposed state compensation linearization method can
avoid the equilibrium point problem of the Jacobian linearization
method and the singularity problem of the feedback linearization method.
Let $y_{\text{d}}=20$, $d=3$, $x_{0}=-1.$ A proportional-integral-differential
(PID) controller is designed for the primary system (\ref{PSys1}),
and the PID parameters are $k_{\text{p}}=0.66,k_{\text{i}}=1.33,k_{\text{d}}=-0.02$,
which is tuned by the ``PID Tuner App'' in the Matlab software.
The simulation result is displayed in Fig. \ref{Fig_Ex1-2}. The output
can track the reference output well. It is noticed that due to continuity,
$y$ goes from $-1$ to $20$, namely $x=0$ at a certain time. This
will make feedback linearization fail. 

\subsection{Example 2: An NMP system with saturation}

\begin{figure}[tbph]
\begin{centering}
\includegraphics[scale=0.9]{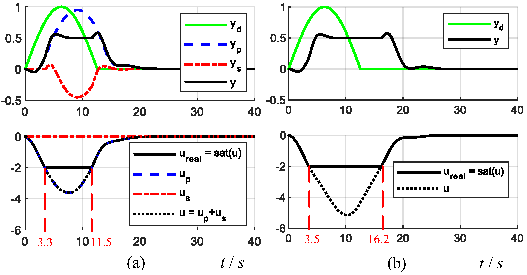} 
\par\end{centering}
\caption{Output response (Example 2): (a) Compensation linearization method;
(b) Jacobian linearization method.}
\label{Fig_Ex4} 
\end{figure}

Consider an NMP system with input saturation as follows 
\begin{align}
\mathbf{\dot{x}} & =\mathbf{Ax}+\mathbf{b}\text{sat}\left(u\right),\mathbf{x}\left(0\right)=\mathbf{x}_{0}\nonumber \\
y & =\mathbf{c}^{\text{T}}\mathbf{x}\label{Ex4}
\end{align}
with 
\[
\mathbf{A=}\left[\begin{array}{ccc}
0 & 1 & 0\\
0 & 0 & 1\\
-4 & -6 & -4
\end{array}\right],\mathbf{b}=\left[\begin{array}{c}
0\\
0\\
1
\end{array}\right],\mathbf{c}=\left[\begin{array}{c}
-1\\
0\\
1
\end{array}\right].
\]
The saturation function sat$\left(\cdot\right)$ is defined as 
\[
\text{sat}\left(u\right)=\left\{ \begin{array}{l}
2,u>2\\
u,2>u>-2\\
-2,u<-2.
\end{array}\right.
\]
The system state is measurable. The control objective is to design
$u$ such that $y$ tracks $y_{\text{d}}$. Next, SCLC is adopted.
System (\ref{Ex4}) can be decomposed into a linear primary system
\begin{align}
\mathbf{\dot{x}}_{\text{p}} & =\mathbf{Ax}_{\text{p}}+\mathbf{b}u_{\text{p}},\mathbf{x}_{\text{p}}\left(0\right)=\mathbf{x}_{0}\nonumber \\
y_{\text{p}} & =\mathbf{c}^{\text{T}}\mathbf{x}_{\text{p}}\label{PSys3}
\end{align}
and a nonlinear secondary system 
\begin{align}
\mathbf{\dot{x}}_{\text{s}} & =\mathbf{Ax}_{\text{s}}+\mathbf{b}\left(\text{sat}\left(u\right)-u\right),\mathbf{x}_{\text{s}}\left(0\right)=\mathbf{0}\nonumber \\
y_{\text{s}} & =\mathbf{c}^{\text{T}}\mathbf{x}_{\text{s}}\label{SSys3}
\end{align}
where let $u_{\text{p}}=u$. According to additive state decomposition,
it follows 
\begin{equation}
{\mathbf{x}}={\mathbf{x}_{\text{p}}}+{\mathbf{x}_{\text{s}},}u=u{_{\text{p}}.}
\end{equation}
Then, the linear system (\ref{PSys3}) ($B=A-C$) can be obtained
by compensating for the nonlinearity of the original system (\ref{Ex4})
($A$) using the secondary system ($C$). If the controller for the
primary system (\ref{PSys3}) is well designed, sat$\left(u_{\text{p}}\right)\rightarrow u_{\text{p}}$,
$u_{\text{p}}\rightarrow u_{\text{d}}$, which also implies $u\rightarrow u_{\text{d}}$.
Then, sat$\left(u_{\text{p}}\right)-u\rightarrow0$. In this case,
since $\mathbf{A}$ is stable, the secondary system (\ref{SSys3})
can stabilize itself. The remaining problem is to design a tracking
controller for (\ref{PSys3}) without saturation. By using state compensation
linearization and the remaining tracking controller for (\ref{PSys3}),
$y\rightarrow y_{\text{d}}$ can hold globally. Let $\mathbf{x}_{0}=\mathbf{0}$
and 
\[
y_{\text{d}}\left(t\right)=\left\{ \begin{array}{cc}
\sin(0.25t), & 0\leq t\leq4\pi\\
0, & 4\pi\leq t.
\end{array}\right.
\]
With a PID controller (the PID parameters are $k_{\text{p}}=-0.5,k_{\text{i}}=-1.3,k_{\text{d}}=0$,
which is tuned by the ``PID Tuner App'' in the Matlab software)
designed for the primary system (\ref{PSys3}), the simulation result
is displayed in Fig. \ref{Fig_Ex4}(a). The simulation result of JLC
is shown in Fig. \ref{Fig_Ex4}(b) for comparison. For the SCLC, the
system undergoes input saturation from $3.3$s to $11.5$s (7.2s duration);
for the JLC, with the same PID controller, the system undergoes input
saturation from $3.5$s to $16.2$s (12.7s duration). It can be seen
that the SCLC can exit saturation earlier than the JLC and can achieve
better tracking performance. In terms of FLC, because the saturation
function is irreversible, the real input $u$ cannot be obtained.
Moreover, because system (\ref{Ex4}) is an NMP system, FLC may lead
to an unstable controller. 

\subsection{Example 3: A mismatched nonlinear system}

Consider a mismatched nonlinear system as follows 
\begin{align}
\dot{x}_{1} & =x_{2}+\sin x_{2}+d_{1},\nonumber \\
\dot{x}_{2} & =-2x_{1}-3x_{2}+2x_{2}^{2}+u+d_{2},\mathbf{x}\left(0\right)=\mathbf{x}_{0}\label{Ex3}
\end{align}
where $\mathbf{x=}\left[x_{1}\ x_{2}\right]^{\text{T}}$is the state,
$u$ is the input, $\mathbf{x}_{0}$ is the initial state value, and
$\mathbf{d=}\left[d_{1}\ d_{2}\right]^{\text{T}}$represents the unknown
uncertainties.

\subsubsection{SCLC}

Based on the proposed state compensation linearization, for the original
system (\ref{Ex3}), the linear primary system (\ref{CL_PSys}) is
chosen as 
\begin{align}
\dot{x}_{\text{p},1} & =2x_{\text{p},2}+d_{1},\nonumber \\
\dot{x}_{\text{p},2} & =-2x_{\text{p},1}-3x_{\text{p},2}+u_{\text{p}}+d_{2},\mathbf{x}_{\text{p}}\left(0\right)=\mathbf{x}_{0}.\label{Ex3p}
\end{align}
Its state-space form is $\mathbf{\dot{x}=}\mathbf{A}\mathbf{x}+\mathbf{B}u+\mathbf{d}$
with 
\[
\mathbf{A}=\left[\begin{array}{cc}
2 & 0\\
-2 & -3
\end{array}\right],\mathbf{B}=\left[\begin{array}{c}
0\\
1
\end{array}\right].
\]
Then, the secondary system (\ref{CL_SSys}) can be obtained as 
\begin{align}
\dot{x}_{\text{s},1} & =2x_{\text{s},2}-x+\text{sin}(x),\nonumber \\
\dot{x}_{\text{s},2} & =-2x_{\text{s},1}-3x_{\text{s},2}+2x{}^{2}+u_{\text{s}},\mathbf{x}_{\text{s}}\left(0\right)=\mathbf{0}.\label{Ex3s}
\end{align}

Next, controller design is carried out. For (\ref{Ex3p}), the LQR
method is employed to calculate a feedback matrix $\mathbf{k}\in\mathbb{R}^{2\times1}$,
and a controller is designed as 
\[
u_{\text{p}}=\mathbf{k}^{\text{T}}\mathbf{\hat{x}}_{\text{p}}.
\]
For (\ref{Ex3s}), a backstepping controller is designed as 
\begin{align*}
u_{\text{s}} & =2\hat{x}_{\text{s},1}+3\hat{x}_{\text{s},2}-2x^{2}\\
 & -a\frac{1}{\hat{x}_{\text{s},1}^{2}+1}\left(\sin\left(\hat{x}_{\text{s},2}\right)+\hat{x}_{\text{s},2}+g\right)-cz_{\text{s},2}
\end{align*}
where $z_{\text{s},2}=\hat{x}_{\text{s},2}+a\arctan\hat{x}_{\text{s},1}$,
$g=\sin\left(x\right)-\sin\left(\hat{x}_{\text{s},2}\right)$, $a,c>0$
are parameters to be specified later. By combining the two controllers,
the final controller follows 
\begin{align*}
u & =u_{\text{p}}+u_{\text{s}}=\mathbf{k}^{\text{T}}\mathbf{\hat{x}}_{\text{p}}+2\hat{x}_{\text{s},1}+3\hat{x}_{\text{s},2}-2x^{2}\\
 & -a\frac{1}{\hat{x}_{\text{s},1}^{2}+1}\left(\sin\left(\hat{x}_{\text{s},2}\right)+\hat{x}_{\text{s},2}+g\right)-cz_{\text{s},2}.
\end{align*}

\subsubsection{JLC for comparison}

By using Jacobian linearization, the following linearized system is
obtained 
\begin{align}
\dot{x}_{1} & =2x_{2}\nonumber \\
\dot{x}_{2} & =-2x_{1}-3x_{2}+u.\label{eq:JL}
\end{align}
Its state-space form is $\mathbf{\dot{x}=}\mathbf{A}\mathbf{x}+\mathbf{B}u$
with 
\[
\mathbf{A}=\left[\begin{array}{cc}
2 & 0\\
-2 & -3
\end{array}\right],\mathbf{B}=\left[\begin{array}{c}
0\\
1
\end{array}\right].
\]
Based on the linear system (\ref{eq:JL}), a controller is designed
as 
\[
u=\mathbf{k}^{\text{T}}\mathbf{x}.
\]

\subsubsection{FLC for comparison}

By using feedback linearization, the following linearized system is
obtained 
\begin{align}
\dot{z}_{1} & =z_{2}\nonumber \\
\dot{z}_{2} & =v\label{eq:FL}
\end{align}
where $z_{1}=x_{1},z_{2}=x_{2}+\sin x_{2}$, $v$ is the virtual input.
Its state-space form is $\mathbf{\dot{z}=}\mathbf{A}\mathbf{z}+\mathbf{B}v$
with 
\[
\mathbf{A}=\left[\begin{array}{cc}
1 & 0\\
0 & 0
\end{array}\right],\mathbf{B}=\left[\begin{array}{c}
0\\
1
\end{array}\right].
\]
Based on the linear system (\ref{eq:FL}), a controller is designed
as 
\[
v=\mathbf{k}^{\text{T}}\mathbf{z}.
\]
Then 
\[
u=\frac{v}{1+\cos x_{2}}+2x_{1}+3x_{2}-2x_{2}^{2}.
\]
A singularity problem may be encountered that the system becomes uncontrollable
and the control input tends to infinity when $x_{2}=-\pi\pm2k\pi,k=1,2,3,$···$,$
i.e., $\cos x_{2}=-1$.

\subsubsection{Robust-Feedback-linearization-based control (RFLC) for comparison}

By using robust feedback linearization, the following linearized system
is obtained 
\begin{align}
\dot{z}_{1} & =2z_{2}\nonumber \\
\dot{z}_{2} & =-2z_{1}-3z_{2}+v\label{eq:RFL}
\end{align}
where $z_{1}=x_{1},z_{2}=0.5x_{2}+0.5\sin x_{2}$. Its state-space
form is $\mathbf{\dot{x}=}\mathbf{A}\mathbf{x}+\mathbf{B}v$ with
\[
\mathbf{A}=\left[\begin{array}{cc}
2 & 0\\
-2 & -3
\end{array}\right],\mathbf{B}=\left[\begin{array}{c}
0\\
1
\end{array}\right].
\]
Based on the linear system (\ref{eq:RFL}), a controller is designed
as 
\[
v=\mathbf{k}^{\text{T}}\mathbf{z}.
\]
Then 
\[
u=\frac{2v}{1+\cos x_{2}}+2x_{1}+3x_{2}-2x_{2}^{2}-\frac{4x_{1}+3x_{2}+3\sin x_{2}}{1+\cos x_{2}}.
\]
Note that, similarly to FLC, a singularity problem may be encountered.

\subsubsection{ADRC for comparison \cite{huang2014active,gao2006active}}

By using ADRC, the following linearized system is obtained 
\begin{align}
\dot{x}_{1} & =x_{2}\nonumber \\
\dot{x}_{2} & =bu\label{eq:ADRC}
\end{align}
where $b$ is the estimation of $1+\cos x_{2}$. Its state-space form
is $\mathbf{\dot{x}=}\mathbf{A}\mathbf{x}+\mathbf{B}u$ with 
\[
\mathbf{A}=\left[\begin{array}{cc}
1 & 0\\
0 & 0
\end{array}\right],\mathbf{B}=\left[\begin{array}{c}
0\\
b
\end{array}\right].
\]
The extended state is 
\[
x_{3}\triangleq\left(1+\cos x_{2}\right)\left(-2x_{1}-3x_{2}+2x_{2}^{2}\right)+\left(1+\cos x_{2}-b\right)u.
\]
An linear extended state observer (LESO) is designed as 
\begin{align*}
\dot{\hat{x}}_{1} & =\hat{x}_{2}-3\omega_{0}\left(\hat{x}_{1}-y\right)\\
\dot{\hat{x}}_{2} & =\hat{x}_{3}-3\omega_{0}^{2}\left(\hat{x}_{1}-y\right)+bu\\
\dot{\hat{x}}_{3} & =-\omega_{0}^{3}\left(\hat{x}_{1}-y\right)
\end{align*}
where the parameter $\omega_{0}$ can be specified later. Then, a
controller is designed as 
\begin{equation}
u=-b^{-1}\hat{x}_{3}+u_{0}\left(\hat{x}_{1},\hat{x}_{2}\right)\label{eq:ADRC-1}
\end{equation}
where $u_{0}\left(\hat{x}_{1},\hat{x}_{2}\right)$ is designed based
on (\ref{eq:ADRC}) as 
\[
u_{0}=\mathbf{k}^{\text{T}}\mathbf{x}.
\]
It can be seen from (\ref{eq:ADRC}) that $b$ must be reversible.
However, in this example, $b$ is the estimation of $1+\cos x_{2}$,
and thus $b$ is irreversible when $x_{2}=-\pi\pm2k\pi,k=1,2,3,$···$,$
i.e., $\cos x_{2}=-1$. It can be concluded that, similarly to FLC
and RFLC, a singularity problem may also be encountered.

\subsubsection{Simulation results}

The LQR method is utilized to determine the feedback matrix $\mathbf{k}$
for all the five control design methods and the parameters are all
selected as $\mathbf{Q}=\text{diag}\left(10,10\right)$, $\mathbf{R}=1$.
Besides, let $a=10$, $c=10$, $\omega_{0}=2$.

The following four scenarios are presented in the simulation to evaluate
and compare the performance of our approach with other linearization
methods.

(i) $\mathbf{x}_{0}=\left[2\ 2\right]^{\text{T}}$, $\mathbf{d=}\mathbf{0}$.

(ii) $\mathbf{x}_{0}=\left[5\ 5\right]^{\text{T}}$, $\mathbf{d=}\mathbf{0}$.

(iii) $\mathbf{x}_{0}=\left[2\ 2\right]^{\text{T}}$, $\mathbf{d=}\left[1\ 1\right]^{\text{T}}$.

(iv) $\mathbf{x}_{0}=\left[2\ 2\right]^{\text{T}}$, $\mathbf{d=}\mathbf{0}$,
a time delay $0.2\text{s}$ is added into the input channel of the
system (\ref{Ex3}).

Two commonly-used performance indices: the integral absolute error
(IAE) and the integral of time-weighted absolute error (ITAE) are
used to evaluate the performance of the above five methods quantitatively.
Please refer to \cite{oliveira2020evolving} to find the computation
method of IAE and ITAE.

\begin{figure}[tbph]
\begin{centering}
\includegraphics[width=0.95\columnwidth]{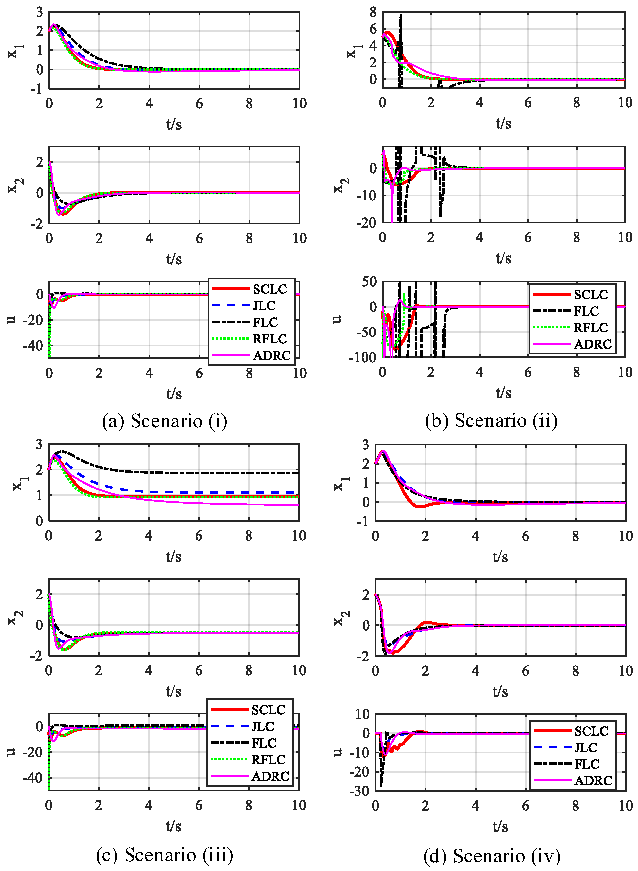} 
\par\end{centering}
\caption{State and input response.}
\label{Fig_Ex2} 
\end{figure}

In Scenario (i), the corresponding state and input response is depicted
in Fig.\ \ref{Fig_Ex2}(a), which shows that all the five methods
can achieve stabilization control. Note that the inputs of FLC and
RFLC have large initial values. From Tab.\ \ref{Tab-1}, it can be
seen that RFLC has the smallest IAE and ITAE values, which is followed
by SCLC in the second place. In Scenario (ii), when the initial state
is increased to $\mathbf{x}_{0}=\left[5,5\right]^{\text{T}}$, the
corresponding state and input response is displayed in Fig.\ \ref{Fig_Ex2}(b).
JLC causes instability and is thus not included in the figure. It
can be found that a singularity phenomenon is incurred by FLC, while
SCLC still works well. RFLC still has large initial input values.
In Scenario (iii), a constant disturbance is considered. From Fig.\ \ref{Fig_Ex2}(c)
and Tab.\ \ref{Tab-1}, one can observe that ADRC has the best performance.
SCLC also performs well. In Scenario (iv), RFLC causes instability
because of the singularity problem and thus do not included in Fig.\ \ref{Fig_Ex2}(d).
Tab.\ \ref{Tab-1} shows that SCLC has the best performance.

It can be concluded that, for system (\ref{Ex3}), JLC has a local
property, and causes instability under large initial state values.
FLC has a singularity problem and does not have enough robustness.
In order to avoid the singularity condition, the working state of
the system must be restricted. RFLC has better robustness than FLC.
ADRC can tackle unknown disturbances better. However, RFLC and ADRC
also have a singularity problem (Some other simulations have also
be done where instability is incurred by the singularity problem in
RFLC and ADRC, especially when the initial values are large). By contrast,
SCLC avoids the local property of JLC, has better robustness than
FLC, and avoids the singularity problem in FLC, RFLC, and ADRC. Besides,
SCLC can accommodate the input time delay better. This is because
all the input time delay is allocated into the primary system by using
state compensation linearization. The secondary system deals with
the known nonlinearity without the effect of the input time delay.
The special strategy mitigates the effect of input time delay on the
system. Thus, state compensation linearization can be complementary
to the existing linearization methods. For some control problems,
if the existing linearization methods do not work well, SCLC is worth
trying. The simulation results verify the theoretical conclusion well.
\begin{table}
\caption{Performance indices of the controllers}

\begin{centering}
\begin{tabular*}{0.95\columnwidth}{@{\extracolsep{\fill}}|@{\extracolsep{\fill}}|c|c|c|c|c|c|c|}
\hline 
Sce.  & Index  & SCLC  & JLC  & FLC  & RFLC  & ADRC\tabularnewline
\hline 
\multirow{2}{*}{(i)} & IAE  & 1.915  & 2.495  & 3.485  & 1.762  & 2.788\tabularnewline
\cline{2-7} \cline{3-7} \cline{4-7} \cline{5-7} \cline{6-7} \cline{7-7} 
 & ITAE  & 1.055  & 1.851  & 3.842  & 0.947  & 4.719\tabularnewline
\hline 
\multirow{2}{*}{(ii)} & IAE  & 5.026  & -  & 12.438  & 3.695  & 5.409\tabularnewline
\cline{2-7} \cline{3-7} \cline{4-7} \cline{5-7} \cline{6-7} \cline{7-7} 
 & ITAE  & 2.866  & -  & 16.546  & 1.960  & 6.840\tabularnewline
\hline 
\multirow{2}{*}{(iii)} & IAE  & 10.908  & 13.129  & 20.027  & 10.485  & 9.626\tabularnewline
\cline{2-7} \cline{3-7} \cline{4-7} \cline{5-7} \cline{6-7} \cline{7-7} 
 & ITAE  & 48.440  & 57.575  & 95.280  & 47.286  & 36.893\tabularnewline
\hline 
\multirow{2}{*}{(iv)} & IAE  & 2.324  & 2.906  & 2.765  & -  & 3.466\tabularnewline
\cline{2-7} \cline{3-7} \cline{4-7} \cline{5-7} \cline{6-7} \cline{7-7} 
 & ITAE  & 1.445  & 2.223  & 2.577  & -  & 6.018\tabularnewline
\hline 
\end{tabular*}
\par\end{centering}
\label{Tab-1} 
\end{table}

\section{Conclusions}

\label{Conclusions}

By virtue of additive state decomposition, this paper proposes a new
linearization process for control design, namely state compensation
linearization. The nonlinear terms can be removed from the original
system through the compensation at the input and state so that a linear
primary system can be obtained. For some control problems like flight
control design with challenging aircraft dynamics, the state compensation
linearization may offer some advantages over traditional Jacobian
linearization and feedback linearization by bridging the gap between
linear control theories and nonlinear systems. Based on the proposed
state compensation linearization, this paper proposes a control framework
for a class of nonlinear systems so that their closed-loop behavior
can follow a specified linear system. Three examples are provided
to show the concrete design procedure. Simulation results have illustrated
that state-compensation-linearization-based control outperforms Jacobian-linearization-based
control and feedback-linearization-based control in some challenging
cases. In future research, state-compensation-linearization-based
stability margin analysis for nonlinear systems would be considered.

\bibliographystyle{IEEEtran}
\bibliography{IEEEabrv}

\end{document}